# Неравенство Йенсена и выпуклые функции
## А.Н. Лепес
РСФМСШИ им.О.А.Жаутыкова, adilsultan01@gmail.com

В книге [1] на странице 218 имеется доказательство неравенства:

$$a^4 + b^4 + c^4 \geq 4abc - 1, \qquad (1)$$

где *a, b, c* – действительные числа. Это доказательство можно записать в одну строчку:

$$a^4 + b^4 + c^4 - 4abc + 1 = (a^2 - 1)^2 + (b^2 - c^2)^2 + 2(bc - a)^2 \geq 0.$$

Интересно только, как догадаться до подобного самостоятельно? Универсального способа нахождения таких представлений, к сожалению, нет. А может можно доказать подобные неравенства как-нибудь по-другому? Ниже представлен еще один способов доказательства неравенства (1).

*Доказательство* (*Лепес А.*). Рассмотрим функцию $f(x) = x^4$, где $x \in R$. Решим предварительно следующую задачу среди функций вида $g(x) = kx^3 + m$ найти такую, что $f(1) = g(1)$, $f'(1) = g'(1)$ и *f*(*x*)≥*g*(*x*), где *x* – любое действительное число. Отсюда, числа *k* и *m* должны удовлетворять следующим равенствам $k + m = 1$, $4 = 3k$. Следовательно, $k = \frac{4}{3}$, $m = -\frac{1}{3}$. Осталось убедиться в том, что *f*(*x*)≥*g*(*x*) при всех *x*∈*R*. Последнее верно, поскольку неравенство $x^4 \geq \frac{4x^3 - 1}{3}$ (см. рис.1) эквивалентно неравенству $(x - 1)^2 \left(3\left(x + \frac{1}{3}\right)^2 + \frac{2}{3}\right) \geq 0$, которое, очевидно верно

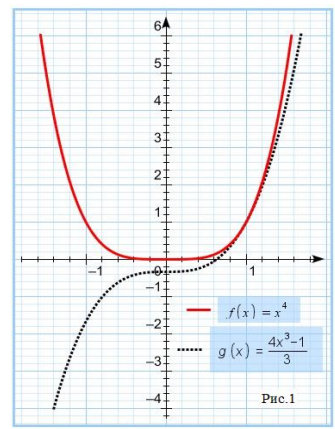

Рис.1

для всякого числа *x*. Отсюда, в силу неравенства Коши $|a|^3 + |b|^3 + |c|^3 \geq 3|abc|$ и неравенства |*t*|≥*t* имеем

$$a^4 + b^4 + c^4 = |a|^4 + |b|^4 + |c|^4 \geq \frac{4(|a|^3 + |b|^3 + |c|^3)}{3} - 1 \geq 4|abc| - 1 \geq 4abc - 1.$$

Таким образом, требуемое неравенство доказано.

Лучше или нет представленное доказательство? С одной стороны, оно существенно длиннее, с другой стороны, оказывается, что подобным способом можно доказывать и другие неравенства, при этом необходимо проводить анализ функции, участвующей в неравенстве, а не изобретать нестандартные преобразования. Например, на Балтийской математической олимпиаде в 2011 году была следующая задача: докажите неравенство

$$\frac{a}{a^3 + 8} + \frac{b}{b^3 + 8} + \frac{c}{c^3 + 8} + \frac{d}{d^3 + 8} \leq \frac{4}{9}, \qquad (2)$$

где *a, b, c, d* – положительные числа, сумма которых равна 4.

На официальной странице Балтийской олимпиады (см.[2]) имеется доказательство неравенства (2), которое вкратце можно записать следующим образом:

$$\frac{a}{a^3+8} + \frac{b}{b^3+8} + \frac{c}{c^3+8} + \frac{d}{d^3+8} =$$

$$= \frac{a}{a^3+1+1+6} + \frac{b}{b^3+1+1+6} + \frac{c}{c^3+1+1+6} + \frac{d}{d^3+1+1+6} \leq [AM-GM] \leq$$

$$\leq \frac{a}{3a+6} + \frac{b}{3b+6} + \frac{c}{3c+6} + \frac{d}{3d+6} = \frac{4 - 2\left(\frac{1}{a+2} + \frac{1}{b+2} + \frac{1}{c+2} + \frac{1}{d+2}\right)}{3} \leq$$

$$\leq [AM-HM] \leq \frac{4 - 2 \cdot \frac{16}{a+2+b+2+c+2+d+2}}{3} = \frac{4}{9}.$$

Здесь мы использовали стандартные обозначения:

AM-GM – неравенство между средним арифметическим и средним геометрическим нескольких чисел;
AM-HM – неравенство между средним арифметическим и средним гармоническим нескольких чисел;

Опять же возникают вопросы: сложно ли самостоятельно догадаться до указанного решения и можно ли доказать это неравенство без искусственных приемов типа 8=1+1+6? В качестве ответа на первый вопрос можно отметить, что Балтийская олимпиада является командным соревнованием, на котором участвуют сборные таких стран как Россия, Германия, Латвия, Финляндия, Швеция, Этония и др. А в 2011 году из 11 команд указанное неравенство смогли доказать только 5 команд. Подобное можно встретить и на других математических олимпиадах. Например, в 2012 году на Международной математической олимпиаде соответствующее неравенство смогли доказать только 189 участников из 548. Основная причина низкого процента учащихся, на наш взгляд, кроется в том, что знание классических неравенств не достаточно для их эффективного применения, в некоторых ситуациях необходимо придумать некоторое преобразование. Поэтому, естественно, возникает желание разработать метод доказательства неравенств хотя бы для некоторого типа, в котором творческая составляющая сведена к минимума. На основании исследований, представленных в настоящем проекте, а также некоторых других разработок, нами разработан новый метод доказательства неравенств: метод отделяющих касательных, описание которого и различные применения приняты в печать в четыре журнала:
1) г. Чебоксары, Россия
Ибатулин И.Ж., Лепес А.Н. Об одном методе доказательства неравенств // Труды XXI конференции «Математика. Образование». 2013. 16 стр. (в печати, октябрь-ноябрь 2013);
2) г. София, Болгария
Ibatulin I.Zh., Lepes A.N. Application of the method of separating tangents to prove inequalities // Didactical Modeling: e-journal 2013. URL: http://www.math.bas.bg/omi/DidMod/index.htm. 13 pages (в печати, 2014);
3) г. Гонконг, Китай

Ibatulin I.Zh., Lepes A.N. Using tangent lines to prove inequalities (part II) // Mathematical Excalibur. 2013-2014. V.1. 6 pages. (в печати, январь-февраль 2014);

4) г. Москва, Россия

Ибатулин И.Ж., Лепес А.Н. Применение касательной для доказательства неравенств // Математика в школе. 2014. 13 стр. (в печати, январь-февраль 2014).

Читателям предлагается попробовать доказать следующие неравенства, тем самым для себя оценить практическую значимость метода отделяющих касательных.

**Пример 1.** ([3], стр. 31, упр. 1.2.9) Пусть даны положительные числа $a, b, c$ такие, что $a^2+b^2+c^2=3$. Докажите неравенство $\frac{1}{a^3+2} + \frac{1}{b^3+2} + \frac{1}{c^3+2} \geq 1$.

**Пример 2.** ([3], стр. 55, упр. 3.1.4) Пусть даны неотрицательные числа $a, b, c, d, e$ такие, что

$$\frac{1}{4+a} + \frac{1}{4+b} + \frac{1}{4+c} + \frac{1}{4+d} + \frac{1}{4+e} = 1.$$

Докажите, что

$$\frac{a}{4+a^2} + \frac{b}{4+b^2} + \frac{c}{4+c^2} + \frac{d}{4+d^2} + \frac{e}{4+e^2} \leq 1.$$

**Пример 3.** (Китай, 2005, [1], стр. 196, задача 132) Пусть даны положительные числа $a, b, c$, сумма которых равна 1. Докажите неравенство

$$10(a^3 + b^3 + c^3) - 9(a^5 + b^5 + c^5) \geq 1.$$

**Пример 4.** ([1], стр. 196, задача 124) Пусть даны положительные числа $a, b, c$ такие, что $a^2+b^2+c^2=12$. Найдите наибольшее значение выражения

$$A = a \cdot \sqrt[3]{b^2 + c^2} + b \cdot \sqrt[3]{c^2 + a^2} + c \cdot \sqrt[3]{a^2 + b^2}.$$

**Пример 5.** (Baltic Way, 2002, задача 4, [4]) Путь дано целое положительное число $n$. Для любых чисел $x_1, x_2, \ldots, x_n \geq 0$ таких, что $\sum_{k=1}^n x_k = 1$, докажите неравенство

$$\sum_{k=1}^n x_k(1-x_k)^2 \leq \left(1 - \frac{1}{n}\right)^2.$$

В качестве еще одной демонстрации метода отделяющих касательных приведем доказательство неравенства (2).

*Доказательство.* Заметим, что если $a=b=c=d=1$, то $\frac{a}{a^3+8} + \frac{b}{b^3+8} + \frac{c}{c^3+8} + \frac{d}{d^3+8} = \frac{4}{9}$. Положим $f(x) = \frac{x}{x^3+8}$, где $x \in (0; 4)$. Составим уравнение касательной к графику функции $f$ в точке $x_0=1$:

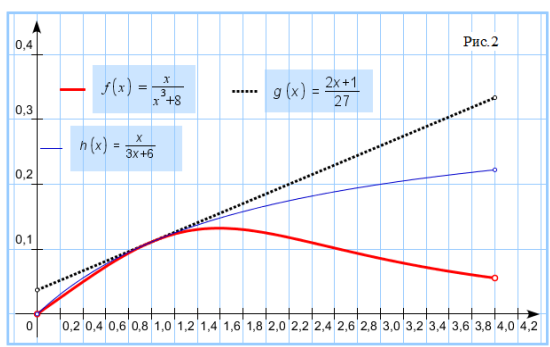

$$y = f(1) + f'(1)(x-1) = \frac{1}{9} + \frac{2}{27}(x-1) = \frac{2x+1}{27}.$$

Покажем, что на интервале (0; 4) график функции *f* лежит не выше касательной прямой $y = \frac{2x+1}{27}$ (см. рис. 2). Действительно, поскольку

$$\frac{x}{x^3+8} \leq \frac{2x+1}{27} \Leftrightarrow (x-1)^2(2x^2+5x+8) \geq 0, \tag{3}$$

то указанное утверждение верно. Применяя неравенство (3), получим

$$\frac{a}{a^3+8} + \frac{b}{b^3+8} + \frac{c}{c^3+8} + \frac{d}{d^3+8} \leq \frac{2(a+b+c+d)+4}{27} = \frac{4}{9}.$$

Тем самым, неравенство (2) доказано.

Справедливости ради, необходимо отметить, что доказательства неравенств с помощью касательной прямой можно найти у некоторых авторов (см., например, [1], [3], [5]).

Важность расположения графика функции по одну сторону от касательной прямой отмечена также в статье В.И. Гаврилова и А.В. Субботина (см.[6]). В статье [6] показана эквивалентность двух определений выпуклых функций: одно через секущую, другое через касательную.

Как известно, для всякой выпуклой на промежутке *I* функции *f* выполнено неравенство Йенсена (см., например, [7]):

$$\frac{\sum_{k=1}^{n} f(x_k)}{n} \geq f\left(\frac{\sum_{k=1}^{n} x_k}{n}\right), \tag{4}$$

где $x_1, x_2, \ldots, x_n \in I$.

Оказывается, что если ослабить требования на функции: рассматривать функции, график которых расположен по одну сторону от касательной прямой, проведённой в заданной точке $x_0$, то мы потеряем выпуклость, однако сохраним неравенство (4), при дополнительном ограничении: $\sum_{k=1}^{n} x_k = nx_0$. В таком случае, будем говорить функция удовлетворяет неравенству Йенсена в заданной точке $x_0$. В более общем случае,

$$\sum_{k=1}^{n} f(x_k) \geq nf(x_0), \tag{5}$$

где $x_1, x_2, \ldots, x_n \in I: \sum_{k=1}^{n} l(x_k) = nl(x_0)$, *l* – некоторая заданная функция.

Между тем, наибольшую сложность на математических олимпиадах представляют доказательства неравенств с невыпуклыми функциями, а в большинстве современной литературы, посвящённой неравенству Йенсена, упоминаются только выпуклые функции. В этом можно убедиться на примере журнала Квант, см. [7]-[9]. Хотя есть и исключения, см., например, [3], стр.68-70. Данная работа призвана также восполнить имеющийся пробел.

Необходимо, отметить, что функция *l* обычно является степенной функцией (см. решение примера 1), но её выбор может быть также обусловлен условием задачи (см. решение пример 2). Возможно, конечно, что график функции лежит по обе стороны от касательной прямой, однако функция удовлетворяет неравенству Йенсена в заданной точке $x_0$ (см. решение примера 3).

*Решение примера* 1 (*Лепес А.*). Заметим, что при $a=b=c=1$ наше неравенство обращается в равенство. Пусть $f(x) = \frac{1}{x^3+2}$, $g(x)=kx^2+m$, где $x \in (0; \sqrt{3})$. Числа $k$ и $m$ выберем такими, что $f(1) = g(1)$, $f'(1) = g'(1)$. То есть $k + m = \frac{1}{3}$, $2k = -\frac{1}{3}$. Следовательно, $g(x) = -\frac{x^2}{6} + \frac{1}{2}$. Неравенство (см. рис. 3) $\frac{1}{x^3+2} \geq -\frac{x^2}{6} + \frac{1}{2}$ является верным, поскольку оно эквивалентно неравенству $x^2(x-1)^2(x+2) \geq 0$. Отсюда,

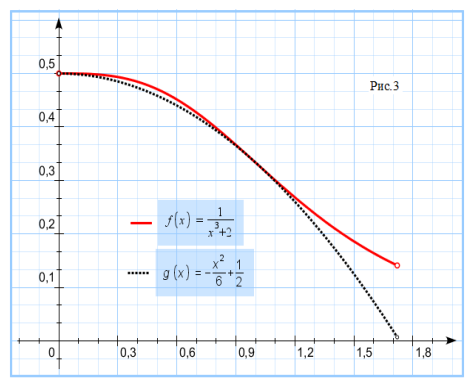

$$\frac{1}{a^3+2} + \frac{1}{b^3+2} + \frac{1}{c^3+2} \geq \frac{3}{2} - \frac{a^2+b^2+c^2}{6} = 1.$$

Таким образом, требуемое неравенство доказано.

*Решение примера* 2 (*Лепес А.*). Положим $f(x) = \frac{x}{4+x^2}$, $g(x) = \frac{k}{4+x} + m$, где $x \geq 0$. Числа $k$ и $m$ выберем таким образом, что $f(1) = g(1)$, $f'(1) = g'(1)$. Отсюда, $k=-3$, $m=0{,}8$. Поскольку неравенство $\frac{x}{4+x^2} \leq \frac{4}{5} - \frac{3}{4+x}$ эквивалентно неравенству $(x-1)^2(x+1) \geq 0$, то оно является верным для всякого $x \geq 0$ (см. рис.4). Отсюда,

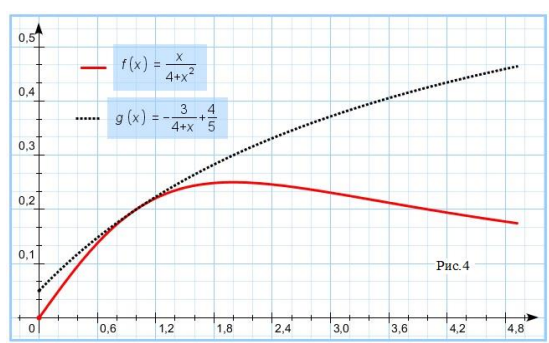

$$\frac{a}{4+a^2} + \frac{b}{4+b^2} + \frac{c}{4+c^2} + \frac{d}{4+d^2} + \frac{e}{4+e^2} \leq 4 - 3\left(\frac{1}{4+a} + \frac{1}{4+b} + \frac{1}{4+c} + \frac{1}{4+d} + \frac{1}{4+e}\right) = 1.$$

Что и требовалось доказать.

*Решение примера* 3 (*Лепес А.*). Пусть дана функция $f(x) = 10x^3 - 9x^5$, $x \in (0; 1]$. Заметим, что при $a = b = c = \frac{1}{3}$ справедливо равенство

$$10(a^3 + b^3 + c^3) - 9(a^5 + b^5 + c^5) = 1.$$

Составим уравнение касательной к графику функции *f* в точке $x_0 = \frac{1}{3}$:

$$y = f\left(\frac{1}{3}\right) + f'\left(\frac{1}{3}\right)\left(x - \frac{1}{3}\right) = \frac{1}{3} + \frac{25}{9}\left(x - \frac{1}{3}\right) = \frac{25}{9}x - \frac{16}{27}.$$

Рассмотрим неравенство (см. рис. 5)

$$10x^3 - 9x^5 \geq \frac{25}{9}x - \frac{16}{27}. \qquad (6)$$

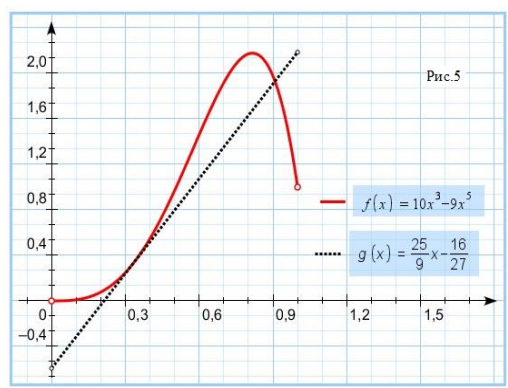

Неравенство (6) можно записать также в следующем виде

$$\left(x - \frac{1}{3}\right)^2 \left(9x^3 + 6x^2 - 7x - \frac{16}{3}\right) \leq 0.$$

А поскольку для всякого $x \in (0; 0{,}9)$ выполнено

$$9x^3 + 6x^2 - 7x - \frac{16}{3} = (x - 0{,}9)(9x^2 + 14{,}1x + 5{,}69) - \frac{637}{3000} < 0,$$

то неравенство (6) справедливо при всех $x \in (0; 0{,}9)$. Значит, если положительные числа *a*, *b*, *c* меньше 0,9 и их сумма равна 1, то согласно неравенству (6) имеем

$$10(a^3 + b^3 + c^3) - 9(a^5 + b^5 + c^5) \geq \frac{25}{9}(a + b + c) - 3 \cdot \frac{16}{27} = 1.$$

Если же хотя бы одно из чисел *a*, *b*, *c* не меньше 0,9, то в силу того, что функция *f* убывает на промежутке [0,9; 1], т.к.

$$f'(x) = -45x^2 \left(x^2 - \frac{2}{3}\right) < 0, \text{ где } x \in [0{,}9; 1],$$

и функции *f* на всем промежутке (0; 1] принимает неотрицательные значения имеем

$$10(a^3 + b^3 + c^3) - 9(a^5 + b^5 + c^5) = f(a) + f(b) + f(c) \geq$$

$$\geq \min_{[0,9;\, 1]} f(x) = f(1) = 1.$$

Тем самым, требуемое неравенство доказано.

Отметим, что решение примера 3 не тривиально. Тем самым, показано, что неравенству Йенсена в заданной точке $x_0$ могут удовлетворять не только выпуклые функции. Если же на это не обращать внимание, то можно совершить ошибку. Так, например, в книге [5] (см., стр. 19) Э.Саппа указал, что функция $f(x) = 10x^3 - 9x^5$ на отрезке [0;1] является выпуклой, и поэтому он утверждал, что она удовлетворяет неравенству Йенсена, т.е.

$$10(a^3 + b^3 + c^3) - 9(a^5 + b^5 + c^5) = f(a) + f(b) + f(c) \geq 3f\left(\frac{1}{3}\right) = 1,$$

где *a*, *b*, *c*>0 и *a*+*b*+*c*=1. Однако, на самом деле, функция *f* является невыпуклой, в этом можно убедиться хотя бы с помощью второй производной: $f''(0{,}1) = 5{,}82$, $f''(0{,}9) = -77{,}22$. Тем самым, функция *f* является примером невыпуклой функции, удовлетворяющей неравенству Йенсена.

Подводя итог, можно отметить, что если необходимо доказать неравенство вида (4), то для этого, оказывается, полезно доказывать вспомогательное неравенство вида

$$f(x) \geq g(x) = kl(x) + m, \tag{7}$$

где числа *k* и *m* выбираются таким образом, чтобы $f(x_0) = g(x_0)$, $f'(x_0) = g'(x_0)$. Выбор точки $x_0$ осуществляется из случая, когда неравенство (4) обращается в равенство. В случае, когда $l(x) = x$, $l(x) = x^2$ или $l(x) = \ln x$ о необходимости указанного способа выбора коэффициентов *k* и *m* указано также в книге Фам Ким Хунга (см. [3], стр.136), однако этот прием работает для любой функции *l*, дифференцируемой в точке $x_0$. Необходимость указанного условия может быть объяснено, например, с помощью теоремы Ферма, поскольку выполнение неравенства (7) в некоторой окрестности точки $x_0$ при условии $f(x_0) = g(x_0)$ означает, что точка $x_0$ является точкой локального минимума функции $h(x) = f(x) - g(x)$, то согласно теореме Ферма имеем

$$h'(x_0) = f'(x_0) - g'(x_0) = 0.$$

Определенный интерес представляет вопрос о применимости касательной прямой не только в случае $l(x)=x$, но и в случаях $l(x) = x^p$ или $l(x) = \ln x$. И здесь могут оказаться полезны некоторые свойства средне степенных (подробнее об этих свойствах см., например, [10], стр. 12-23). Начнем с определения.

Пусть даны положительные числа $x_1, x_2, \ldots, x_n$. Если число $\alpha$ отлично от нуля, то средним степенным чисел $x_1, x_2, \ldots, x_n$ называется число

$$c_\alpha(x_1, x_2, \ldots, x_n) = \left(\frac{\sum_{j=1}^n x_j^\alpha}{n}\right)^{\frac{1}{\alpha}}.$$

Если $\alpha=0$, то

$$c_0(x_1, x_2, \ldots, x_n) = \sqrt[n]{x_1 \cdot x_2 \cdot \ldots \cdot x_n}.$$

Справедлива

**Теорема 1 (достаточные условия для выполнения неравенство Йенсена).** *Пусть даны число $\alpha$, положительное число $x_0$ и натуральное число $n \geq 2$. Пусть дана функция $f$, определенная на множестве всех положительных чисел, и дифференцируема в точке $x_0$. Если $(\alpha - 1) \cdot f'(x_0) \leq 0$ и для всякого положительного числа $x$ такого, что $x^\alpha < nx_0^\alpha$ выполнено $f(x) \geq f(x_0) + f'(x_0)(x - x_0)$, то*

$$\sum_{j=1}^n f(x_j) \geq nf(x_0),$$

*где $x_1, x_2, \ldots, x_n > 0$ и $c_\alpha(x_1, x_2, \ldots, x_n) = x_0$.*

*Доказательство.* Если $\alpha \neq 0$ и для положительных чисел $x_1, x_2, \ldots, x_n > 0$ выполнено $c_\alpha(x_1, x_2, \ldots, x_n) = x_0$, то вне зависимости положительно или нет число $\alpha$ для каждого $k=1, 2, \ldots, n$ выполнено неравенство $x_k^\alpha < nx_0^\alpha$. Если же $\alpha=0$, то указанное неравенство также выполнено.

Пусть $\alpha<1$. Тогда $f'(x_0) \geq 0$. Согласно теореме о монотонности средне степенных (см., например, [10], стр. 21) имеем

$$\sum_{j=1}^n f(x_j) \geq \sum_{j=1}^n \left(f(x_0) + f'(x_0)(x_j - x_0)\right) = nf(x_0) + nf'(x_0)(c_1(x_1, x_2, \ldots, x_n) - x_0) \geq$$

$$\geq nf(x_0) + nf'(x_0)(c_\alpha(x_1, x_2, \ldots, x_n) - x_0) = nf(x_0)$$

Если $\alpha>1$, то $f'(x_0) \leq 0$ и $nf'(x_0)c_1(x_1, x_2, \ldots, x_n) \geq nf'(x_0)c_\alpha(x_1, x_2, \ldots, x_n)$. Отсюда, также как в предыдущем случае следует неравенство $\sum_{j=1}^n f(x_j) \geq nf(x_0)$.

В случае $\alpha=1$ доказательство тривиально:

$$\sum_{j=1}^{n} f(x_j) \geq \sum_{j=1}^{n} \left( f(x_0) + f'(x_0)(x_j - x_0) \right) = n f(x_0).$$

Тем самым, теорема 1 доказана.

*Замечание* 1. В случае, если α=1 в теореме 1 утверждается, что для всякой выпуклой вниз функции *f* выполнено неравенство Йенсена в точке $x_0$. Это общеизвестный факт, который упоминается, например, в статье [7].

Из теоремы 1 следует важный практический смысл. Если в условии задачи имеется ограничение, например, на сумму квадратов переменных, то рассматривать касательную прямую имеет смысл, если производная соответствующей функции в точке касания отрицательна. Если же, например, произведение переменных постоянно, то желательно, чтобы производная функции в точке касания была положительной. Если же указанное условие для функции не выполнено, то желательно в качестве касательной брать многочлен той же степени, что и порядок степени переменных в условии задачи.

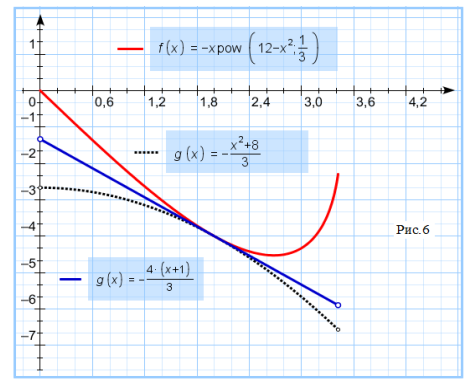

Рис.6

*Решение примера* 4 (*Лепес А.*). Пусть $f(x) = -x \cdot \sqrt[3]{12 - x^2}$, где $x \in (0; 2\sqrt{3})$. Составим уравнение касательной к графику функции *f* в точке $x_0$=2:

$$y = f(x_0) + f'(x_0)(x - x_0) = -4 - \frac{2}{3}(x-2) = -\frac{4x+4}{3}.$$

Заметим, что для всякого положительного числа $x$ выполнено неравенство $f''(x) = \frac{\frac{2}{3}x}{\sqrt[3]{(12-x^2)^2}} + \frac{\frac{16}{9}x^3}{\sqrt[3]{(12-x^2)^5}} > 0$, следовательно, функция *f* выпукла вниз на интервале $(0; 2\sqrt{3})$, а значит ее график лежит не ниже касательной прямой (см. рис.6). Тогда согласно теореме 3 при α=2, $x_0$=2, *f*(2)=-4, *f* '(2)=-2/3 имеем

$$-a \cdot \sqrt[3]{b^2 + c^2} - b \cdot \sqrt[3]{c^2 + a^2} - c \cdot \sqrt[3]{a^2 + b^2} \geq -3 \cdot 4 = 12,$$

причем равенство достигается при *a*=*b*=*c*=2. Таким образом,

$$\max_{\substack{a^2+b^2+c^2=3, \\ a,b,c>0}} \left( a \cdot \sqrt[3]{b^2 + c^2} + b \cdot \sqrt[3]{c^2 + a^2} + c \cdot \sqrt[3]{a^2 + b^2} \right) = 12.$$

*Замечание* 2. В решении примера 4 мы использовали касательную прямую, однако можно было использовать и параболу (см.рис.6). То есть был выбор, а в решении примера 1 мы не рассматривали касательную прямую, однако в точке касания производная функции – отрицательна, что соответствует условию теоремы 1. Единственное, что из рисунка 3 можно предположить, что график функции не лежит по одну сторону от касательной прямой, поэтому необходимо разбить область определения также как в решении примера 3. Ясно, что такое доказательство, даже если оно возможно, будет сложнее представленного выше.

Интересно, что теорема 1 может подсказать не только, когда стоит рассматривать касательную прямую, но и когда необходимо выбрать что-то другое. Например, в 2005 году на Балтийской олимпиаде было следующее неравенство (задача 5):

$$\frac{a}{a^2+2} + \frac{b}{b^2+2} + \frac{c}{c^2+2} \leq 1,$$

где *a, b, c* – положительные числа такие, что *abc*=1.

Из 11 команд данное неравенство смогли доказать только команды стран Эстонии, Исландии и Польши. При этом в книге [1] на 246 (задача 56) имеется решение, в котором используется только неравенство Коши. Однако опять же для применения неравенства Коши необходимы представления типа $a^2+2=a^2+1+1$, которые, безусловно, понятны, не ясно только как догадаться до них самостоятельно. Ведь 8 команд не смогли доказать указанное неравенство. Благодаря же теореме 1, хотя она и не является критерием, можно обратить внимание на то, что производная функции $f(x) = -\frac{x}{x^2+2}$, где *x*>0, в точке касания $x_0$=1 должна быть положительна, поскольку по условию α=0. А эта производная отрицательна. Следовательно, рассматривать в качестве локальной опорной кривой прямую, несмотря на то, что функция *f* выпукла вниз, нет необходимости, значит, желательно попробовать другой тип кривой. И как было указано для подобных случаев у Фам Ким Хунга в книге [3] на странице 136, можно попробовать следующую функцию

$$g(x) = -\frac{3 + lnx}{9}, \text{где } x > 0.$$

Конечно, график функции *f* не лежит по одну сторону от графика функции *g*, однако поступая также, как и в решении примера 3, можно доказать требуемое неравенство.

В теореме 1 ограничение имеется на сумму некоторых степеней переменных, а в качестве касательная рассматривается прямая. Естественно поставить аналогичный вопрос о целесообразности рассмотрения в качестве касательной некоторую степенную функцию, если ограничение в условии задачи имеется на сумму переменных.

**Теорема 2 (достаточные условия для выполнения неравенство Йенсена).** *Пусть даны число α≠0, положительное число $x_0$ и натуральное число n≥2. Пусть дана функция f, определенная на множестве всех положительных чисел, и дифференцируема в точке $x_0$. Если $(\alpha - 1) \cdot f'(x_0) \geq 0$ и для всякого положительного числа x такого, что $x < nx_0$ выполнено $f(x) \geq \frac{f'(x_0)}{\alpha \cdot x_0^{\alpha-1}}(x^\alpha - x_0^\alpha) + f(x_0)$, то*

$$\sum_{j=1}^{n} f(x_j) \geq nf(x_0),$$

*где $x_1, x_2, ..., x_n$>0 и $x_1 + x_2 + \cdots + x_n = nx_0$.*

*Доказательство* теоремы 2 аналогично доказательству теоремы 1.

В заключении остановимся на одном частном случае.

**Теорема 3.** *Пусть дан многочлен P(x)=$ax^3+bx^2+cx+d$, a≠0, натуральное число n и положительное число $x_0$. Пусть даны неотрицательные числа $x_1, x_2, …, x_n$, сумма которых равна $nx_0$. Если $2ax_0+b≥0$ и $(n+2)ax_0+b≥0$, то справедливо неравенство Йенсена в точке $x_0$.*

*Доказательство.* Заметим, что

$$\sum_{j=1}^n P'(x_0)(x_j - x_0) = P'(x_0)\left(\sum_{j=1}^n x_j - nx_0\right) = 0,$$

отсюда, неравенство $\sum_{j=1}^n P(x_j) \geq nP(x_0) = nP(x_0) + \sum_{j=1}^n P'(x_0)(x_j - x_0)$ эквивалентно следующему

$$\sum_{j=1}^n \left(P(x_j) - P(x_0) - P'(x_0)(x_j - x_0)\right) = \sum_{j=1}^n (x_j - x_0)^2 (ax_j + 2ax_0 + b) \geq 0.$$

Заметим, что из того, что сумма неотрицательных чисел $x_1, x_2, \ldots, x_n$ равна $nx_0$ следует, что каждое из этих чисел принадлежит промежутку $[0; nx_0]$. Однако в силу монотонности функции $ax + 2ax_0 + b$ и неотрицательности чисел $2ax_0+b$ и $(n+2)ax_0+b$ для всякого $x \in [0; nx_0]$ справедливо неравенство $ax + 2ax_0 + b \geq 0$. Последнее означает, что требуемое неравенство является верным. Следовательно, теорема 3 доказана.

*Замечание* 3. Как известно, для того, чтобы многочлен $P(x)=ax^3+bx^2+cx+d$ был выпуклым вниз необходимо и достаточно, чтобы $P''(x) = 6ax + 2b \geq 0$ на всем промежутке $[0; nx_0]$. А для этого необходимо, чтобы $b\geq0$ и $3nax_0+b\geq0$. Если числа $a, b, n, x_0$ удовлетворяют указанным неравенствам, то они и удовлетворяют неравенствам $2ax_0+b\geq0$ и $(n+2)ax_0+b\geq0$. Что вполне соответствует общеизвестному факту, что неравенство Йенсена выполнено для выпуклых функций. Однако, если числа $a, b, n, x_0$ удовлетворяют неравенствам $2ax_0+b\geq0$ и $(n+2)ax_0+b\geq0$, то они не обязаны удовлетворять обоим неравенствам $b\geq0$ и $3nax_0+b\geq0$. Например, такими числами могут быть $n=3, a=1, b=-1, x_0=1$. То есть любая функция вида $P(x) = x^3 - x^2 + cx + d$ является невыпуклой на отрезке $[0; 3]$, но при этом для нее справедливо неравенство Йенсена в точке $x_0=1$. Ясно, что такой набор чисел $n, a, b, x_0$ не единственный, можно, например, взять $a=1, b= -x_0, x_0>0, n\geq2$. Тем самым, показано, что существует бесконечно много невыпуклых функций, для которых выполнено неравенство Йенсена в заданной точке.

*Решение примера* 5 (*Лепес А.*). В случае $n\geq2$ требуемое неравенство следует из теоремы 5 при $a$=-1, $b$=2, $c$=-1, $d$=0, $x_0 = \frac{1}{n}$. Если $n$=1, то наше неравенство запишется в виде 0≤0, что является верным.

## Литература